\input amstex
\documentstyle{amsppt}
\loadbold
\pageheight{7.0in}

\magnification=\magstep 1
\CenteredTagsOnSplits
%
%  Macros required because they aren't yet in amsppt.sty; see
%  Guidelines for Preparing Electronic Manuscripts (AMS-TeX).
%  Macro for current address.
\def\curraddr#1\endcurraddr{\address {\it Current address\/}: #1\endaddress}

%The next three lines get rid of the amstex logo at the
%bottom of page 1, by redefining it as a null string.
%It still works if you comment out the third line, for some
%reason. See Texbook p. 37 for a listing of category codes.
\catcode`\@=11
\redefine\logo@{}
\catcode`\@=13

\def\<{\left<}
		%for inner products
\def\>{\right>}

\topmatter
\title
Extremal marginal tracial states in coupled systems
\endtitle

\author Geoffrey L. Price$^{\ast}$ and Sh\^oichir\^o Sakai
\endauthor
\affil
Department of Mathematics 9E\\
United States Naval Academy \\
Annapolis, MD 21402, USA\\\\
5-1-6-205, Odawara \\
Aoba-Ku \\
Sendai, Japan 980-0003
\endaffil
\thanks $^{\ast}$Research supported by NSF grant DMS-0400841. 
\endthanks

\abstract

Let $\Gamma$ be the convex set consisting of all
states $\phi$ on the tensor product $B\otimes B$ of 
the algebra $B = M_n(\Bbb C)$ of all $n\times n$ 
matrices over the complex numbers $\Bbb C$ with the 
property that 
the restrictions
$\phi_{\restriction{B\otimes I}}$ and
$\phi_{\restriction{I\otimes B}}$ are
the unique tracial states on $B\otimes I$ and
$I\otimes B$.  We find necessary and sufficient 
conditions for such a state, called a marginal tracial
state, to be extremal in $\Gamma$.  We also give a
characterization of those extreme points in $\Gamma$
which are pure states.
We conjecture that
all extremal marginal tracial states are pure states.
\endabstract
\endtopmatter
\document

\subhead
{1. Introduction}
\endsubhead

For $n\geq 2$ let $B$ be the type 
$I_n$-factor consisting of the 
$n\times n$ matrices $M_n(\Bbb C)$ over the complex numbers 
$\Bbb C$. Let $B\otimes B$ be the tensor product of $B$ with 
itself.  Then $B\otimes B$ is isomorphic to the full matrix 
algebra of $n^2\times n^2$ matrices over $\Bbb C$.  
The mapping $x\rightarrow x\otimes I$ (respectively, $x\rightarrow 
I\otimes x$), for 
$x\in B$ gives a unital embedding of $B$ into the subalgebra 
$B\otimes I$ (respectively, $I\otimes B$) of $B\otimes B$.
 
Let $\tau$ be the unique tracial state on $B\otimes B$.  
By restriction to $B\otimes I$ (respectively, to $I\otimes B$), 
$\tau(x\otimes I) = \tau_B(x)$ (respectively, $\tau(I\otimes x) = 
\tau_B(x)$).  Under these identifications we will abuse notation 
somewhat by using $\tau$ interchangeably to refer to the tracial 
state on $B\otimes B$ as well as on $B$.

A state $\phi$ on $B\otimes B$ is called a marginal tracial state 
if the restrictions of $\phi$ to $B\otimes I$ and $I\otimes B$ 
are the tracial states on $B\otimes I$ and $I\otimes B$ 
respectively.  Note that the set $\Gamma$ consisting of all 
marginal tracial states on $B\otimes B$ is non-empty (since, of 
course, $\tau\in \Gamma$) and convex.  Moreover, $\Gamma$ is 
$\sigma(T\otimes T,B\otimes B)$-compact, where $T$ is the dual of 
$B$.

In \cite{P1} K. R. Parthasarathy has shown for 
the case $n=2$ that any extremal marginal tracial state is 
a pure state.
Our work in the present paper is directed towards 
determining whether or not 
one can extend Parthasarathy's result to all $n$.  Although we 
cannot solve the problem we have obtained some partial results 
which we hope will prove useful.  
Our results suggest that 
Parthasarathy's result may indeed extend to all $n$ (see
Lemma $2.7$ and the remark following).

Given these considerations it is obviously of interest
to determine those pure states on $B\otimes B$ which are
marginal tracial states.  The proof
of the theorem below shows that there
is a one-to-one correspondence between such states and 
$SU(n)$.  The proof, which we include
for the sake of completeness, is achieved by constructing the 
Schmidt decomposition of a unit vector corresponding to
the pure state, \cite{Sc}, \cite{VW}, \cite{EK}.

\proclaim{Theorem 1.1}  A marginal tracial state $\phi$ on $B\otimes 
B$ is a pure state if and only if there are orthonormal bases
 $\{ f_i: 1\leq i \leq n\}$ and $\{g_i:1\leq i \leq n \}$ of $\Bbb C^n$ 
such that
$\phi(A) = \left<A\xi,\xi\right>$ for 
all $A\in B\otimes B$, where $\xi= \sum_{i=1}^n 
\frac{f_i\otimes g_i}{\sqrt{n}}$.
\endproclaim

\remark{Remark 1.1}  A decomposition of a vector $\eta$ in $\Bbb 
C^n\otimes \Bbb C^n$ into the form $\sum_{i=1}^n a_if_i\otimes 
g_i$, where $\{f_i:i=1,\dots,n\},\{g_j:j=1,\dots,n\}$
are a pair of orthonormal bases of $\Bbb C^n$, 
 is known as a Schmidt decomposition, (see \cite{EK} for the 
reference to \cite{Sc})).  In \cite{EK} it is shown that 
for any pure state $\omega$ 
on $M_p(\Bbb C)\otimes M_q(\Bbb C)$ there is an orthonormal
basis $\{f_i: i=1,\dots,p\}$ for $\Bbb C^p$ and another,
$\{g_j:j=1,\dots,q\}$  for $\Bbb C^q$ such that
the unit vector 
$\eta$ satisfying $\omega = \left<\cdot \eta,\eta \right>$ 
admits a Schmidt decomposition $\eta = \sum_{i=1}^r a_if_i\otimes 
g_i$, where $r \leq \text{min}(p,q)$.

In the terminology of quantum computing the vector 
$\xi$ in the statement of the theorem is said to be a maximally 
entangled vector which means that $r$ above is maximal
(i.e., $r=n$) in the Schmidt decomposition of $\xi$.
\endremark

\demo{Proof}  If $\{f_i\}$ and $\{g_i\}$ are orthonormal bases 
for $\Bbb C^n$ it is straightforward to verify that the pure 
state $(\cdot \xi, \xi)$ is a marginal tracial state on $B\otimes B$, 
where $\xi = \sum_{i=1}^n \frac{f_i\otimes 
g_i}{\sqrt{n}}$.  Conversely,
suppose $\phi$ is a  
marginal tracial state
that is also pure.  
Then for $1\leq i,j \leq n$ there are complex numbers
 $\lambda_{ij}\in \Bbb 
C$ be such that $\phi(A) = \left< A\xi,\xi \right>$ 
where $\xi = \sum_{i,j=1}^n 
\lambda_{ij}f_i\otimes f_j$.  For $x\in B$ we have 
$\phi(x\otimes I) = \left<(x\otimes I)\xi,\xi \right> 
= \tau(x\otimes I)$.  Let 
$\{e_{rs}: 1\leq r,s \leq n\}$ be matrix units of $B$ 
corresponding to the orthonormal basis $\{f_i\}$, i.e. 
$e_{rs}f_i= \delta_{is}f_r$, then clearly
$$
\delta_{rs}/n=\tau(e_{rs}) = 
\phi(e_{rs}\otimes I) = \left<(e_{rs}\otimes I)\xi,\xi \right> 
= \sum_{j=1}^n 
\lambda_{sj}\bar \lambda_{rj},
$$
so that $\Lambda = \left(\sqrt{n}\lambda_{ij} \right)$ is 
a unitary 
matrix.  
For $1\leq i \leq n$ set $g_i = \sum_{j=1}^n \sqrt{n}\lambda_{ij}f_j$.  
Since $\Lambda$ is a unitary matrix the set $\{g_i: 1\leq i \leq 
n\}$ is an orthonormal basis of $\Bbb C^n$, and $\phi = 
\left< \cdot \xi,\xi \right>$ where $\xi = \sum_{i=1}^n 
\frac{f_i\otimes 
g_i}{\sqrt{n}}$. \qed
\enddemo
\remark{Remark $1.2$}  From the proof above it 
follows that if $\phi$ is a pure state on $B\otimes B$ that 
restricts to the tracial state on $B\otimes I$, 
i.e. $\phi(x\otimes I)= \tau(x)= \tau_B(x)$ for all $x\in B$, then 
$\phi$ is automatically a marginal tracial state
on $B\otimes B$.  For  
it follows under these hypotheses that $\phi = \left< \cdot 
\xi,\xi \right>$ where $\xi = \sum_{i=1}^n 
\frac{f_i\otimes 
g_i}{\sqrt{n}}$ and so for any $y\in B$,  $\phi(I\otimes 
y) = \left<(I\otimes y)\xi,\xi \right> = 
\left<(I\otimes y)(\sum_{i=1}^n 
\frac{f_i\otimes 
g_i}{\sqrt{n}}), \sum_{i=1}^n 
\frac{f_i\otimes 
g_i}{\sqrt{n}}\right> = 
\sum_{i=1}^n \left< yg_i,g_i \right>/n = \tau(y)$.
\endremark

\remark{Remark $1.3$} Using an argument similar to
the proof above (see also 
\cite{EK}) one can show that if $\phi = 
\left< \cdot \xi, \xi \right>$ is a 
pure
marginal tracial state then the vector $\xi$ has the following
property:  let $\eta$ be any vector in $\Bbb C^n 
\otimes \Bbb C^n = \Cal H\otimes \Cal H$, the vector space on which 
$B\otimes B$ acts; then there is a linear operator $\Phi\in B$ 
such that $(\Phi\otimes I)\xi = \eta$.
\endremark

\subhead{2. Main Results}
\endsubhead
Let $\frak S$ be the convex set of all states on $B\otimes B$; 
then for each $\phi\in \frak S$ there is a unique positive 
element $h_{\phi}\in T\otimes T$ such that $\phi(a)= 
\tau(ah_{\phi})$, all $a\in B\otimes B$, and 
$\tau(h_{\phi})=\phi(I)=1$.  Conversely, let $h$ be any 
element of $T\otimes T$ satisfying $\tau(h)=1=\Vert h \Vert_1$, 
where $\Vert \cdot \Vert_1$ denotes the trace norm on $T\otimes 
T$; then $h$ must be positive and 
$\phi^h(a)=\tau(ah), a\in B\otimes B$, is a state on 
$B\otimes B$, \cite{Sa}.  
Under this correspondence we shall identify the 
set of all states on $B\otimes B$ with the set of all positive 
elements of $B\otimes B$ with trace norm $1$. 

Let $P$ (resp., Q) be the canonical conditional expectation of 
$B\otimes B$ onto $B\otimes I$ (resp., $I\otimes B$); then the 
linear mappings $P$ and $Q$ satisfy the following conditions:
\roster
\item"{(i)}" $P(h)\geq 0$ (resp. $Q(h)\geq 0$) for $h\geq 0$ in 
$B\otimes B$, and $P(I)=I (resp. Q(I)=I)$.
\item"{(ii)}" $P(axb)=aP(x)b$ for $a,b\in B\otimes I$ and $x\in 
B\otimes B$, (resp., $Q(cxd)=cQ(x)d$ for $c,d\in I\otimes B$ and 
$x\in B\otimes B$).
\item"{(iii)}" $P(xk)=P(kx)$ for $k\in I\otimes B$ and $x\in 
B\otimes B$ (resp. $Q(x\ell) = Q(\ell x)$ for $\ell\in B\otimes 
I$ and $x\in B\otimes B$).
\item"{(iv)}" $P(x^*x)=0$ (resp. $Q(x^*x)=0$) if and only if 
$x=0$, for $x\in B\otimes B$.
\item"{(v)}" $\tau(P(x)) = \tau(Q(x)) = \tau(x)$ for $x\in 
B\otimes B$.
\item"{(vi)}" $P(x^*)=P(x)^*$ and $Q(x^*)=Q(x)^*$ for $x\in 
B\otimes B$.
\item"{(vii)}" $\Vert P(x) \Vert \leq \Vert x \Vert$ and 
$\Vert Q(x) \Vert \leq \Vert x \Vert$ for all $x\in B\otimes B$.
\endroster

Since $T\otimes T$ and $B\otimes B$ are identical as sets, the 
conditional expectation $P$ (resp. $Q$) may be viewed as a 
positive linear mapping of $T\otimes T$ onto $T\otimes I$ (resp. 
of $T\otimes T$ onto $I\otimes T$) and satisfying the same 
properties above.  Concerning the norm we have 
$\Vert P(x) \Vert_1 \leq \Vert x \Vert_1$ and 
$\Vert Q(x) \Vert_1 \leq \Vert x \Vert_1$ for $x\in T\otimes T.$  
In fact,
$$
\align
\Vert P(x) \Vert_1 
 &= \underset{\Vert a \Vert \leq 1,a\in B\otimes 
I}\to{\text{sup}}
\vert \tau (aP(x))\vert \\
 &= \underset{\Vert a \Vert \leq 1,a\in B\otimes 
I}\to{\text{sup}} 
\vert \tau (P(aP(x)))\vert \\ 
&= \underset{\Vert a \Vert \leq 1,a\in B\otimes I}\to{\text{sup}}
\vert \tau (P(a)P(x))\vert   \\
&= \underset{\Vert a \Vert \leq 1,a\in B\otimes I}\to{\text{sup}} 
\vert \tau (P(a)x)\vert \\ 
& = \underset{\Vert a \Vert \leq 1,a\in B\otimes 
I}\to{\text{sup}} 
\vert \tau (ax)\vert      \\
& \leq 
\underset{\Vert b \Vert \leq 1,b\in B\otimes 
B}\to{\text{sup}} \vert \tau (bx)\vert  \\
& \leq \Vert x \Vert_1.  \\
\endalign
$$

\proclaim{Lemma 2.1}  A state $h$ on $B\otimes B$ is a marginal 
tracial state if and only if $P(h)=Q(h)=I$.
\endproclaim
\demo{Proof}  Suppose $h$ is a marginal tracial state; then 
$\tau(ah) = \tau(a)$, $a\in B\otimes I$ is the tracial state on 
$B\otimes I$ and $\tau(ah)= \tau(P(ah)) = \tau(aP(h)) = 
\tau(a)$, $a\in B\otimes I$.  Hence $P(h) = I$, and 
analogously, $Q(h)=I$.  Conversely, if $P(h)=Q(h)=I$ for a 
positive element $h$ in $T\otimes T$ then 
$\tau(ah)=\tau(P(ah))=\tau(aP(h))=\tau(a)$ for $a\in B\otimes 
I$ and so $h$ is tracial on $B\otimes I$, and analogously it is 
tracial on $I\otimes B$, so that $h$ is a marginal tracial 
state.\qed
\enddemo

Consider the linear mapping $P+Q$ on $T\otimes T$ and let 
$V=ker(P+Q)$; then for $v\in V$, $(P+Q)(v) = 0$ and 
$P(P+Q)(v) = P(v)+PQ(v) = P(v)+\tau(v)I=0$.  Analogously, 
$Q(v)+\tau(v)I=0$; hence $(P+Q)(v) +2\tau(v)I=0$; therefore 
$\tau(v)=0$ and so $P(v)=Q(v)=0$.  

\remark{Remark 2.1} It is straightforward to see 
that $V = N\otimes N$ where $N\subset B$ is the
vector space of all elements $b\in B$ with trace $0$.
\endremark

\proclaim{Lemma 2.2}  $(I-V)\cap (T\otimes T)^+$ is the set of 
all marginal tracial states on $B\otimes B$, where $(T\otimes 
T)^+$ is the set of all positive elements of $T\otimes T$.
\endproclaim
\demo{Proof}  Let $h\in (I-V)\cap (T\otimes T)^+$ and put $h=I-
v$, $v\in V$.  Since $h\geq 0, \Vert h \Vert_1 = \tau(h) = 
\tau(I-v)=\tau(I)=1$; hence $h$ is a state on $B\otimes B$.  
Moreover, $P(I-v)= P(I)-P(v)=I$ and $Q(I-v)=I$ and therefore,
by the previous lemma, $h$ is a 
marginal tracial state.

Conversely suppose $h$ is a marginal tracial state on $B\otimes 
B$; then $P(I-h)= I-I =0$ and $Q(I-h)=0$; hence $I-h\in V$ and 
$h= I-(I-h)$.\qed
\enddemo

\proclaim{Theorem 2.3}  Let $h_0$ be a marginal tracial state on 
$B\otimes B$.  Then the following conditions are equivalent:
\roster
\item"{(i)}"  $h_0$ is extremal among the 
marginal tracial states on $B\otimes B$,
\item"{(ii)}" $(R(h_0)T\otimes TR(h_0))\cap V = \{0\},$
\item"{(iii)}" $(R(h_0)T\otimes TR(h_0))\cap \{\lambda I+V: 
\lambda \in \Bbb C\} = \{\lambda h_0: \lambda \in \Bbb C\}$,
\endroster
where $R(h_0)$ is the range projection of $h_0$.
\endproclaim
\demo{Proof} $(i)\implies (ii)$:  Suppose there is a nonzero 
element $v$ in $(R(h_0)T\otimes TR(h_0))\cap V$.  Since 
$R(h_0)T\otimes TR(h_0)$ and $V$ are selfadjoint, $v^*\in
(R(h_0)T\otimes TR(h_0))\cap V$ and so $v+v^*, iv-
iv^*\in (R(h_0)T\otimes TR(h_0))\cap V$; hence without loss of 
generality we may assume that $v$ is self-adjoint.  Then there is 
a positive real number $\lambda$ such that
$
-\lambda R(h_0) \leq v \leq \lambda R(h_0).$  Therefore there 
exists a positive number $\mu$ such that $-\mu h_0 \leq -\lambda 
R(h_0) \leq v \leq \lambda R(h_0) \leq \mu h_0$.  Hence 
$\mu h_0 \pm v \geq 0$ and so $h_0 \pm \frac{v}{\mu} \geq 0$.  
Then $\Vert h_0 \pm \frac{v}{\mu} \Vert_1 = \tau(h_0 \pm 
\frac{v}{\mu}) = \tau(h_0) = 1$.  Moreover, $P(h_0 \pm 
\frac{v}{\mu})= P(h_0)=I$                   
(resp. $Q(h_0 \pm 
\frac{v}{\mu})= Q(h_0)=I)
$and so $h_0 \pm \frac{v}{\mu}$ are 
marginal tracial states and $h_0 = \frac{(h_0+\frac{v}{\mu})+(h_0-
\frac{v}{\mu})}{2}$, a contradiction.

$(ii)\implies (iii)$:  Clearly  
$h_0\in R(h_0)T\otimes TR(h_0)$ and 
$h_0 = I + (h_0-I)$ with $h_0-I\in V$; hence $\{ \lambda h_0: 
\lambda \in \Bbb C\} \subset
(R(h_0)T\otimes TR(h_0))\cap \{\lambda I+V: 
\lambda \in \Bbb C\}$.  If there is an element $k$ in
$(R(h_0)T\otimes TR(h_0))\cap \{\lambda I+V: 
\lambda \in \Bbb C\}$ that is not a scalar multiple of $h_0$ then 
$k= \lambda_1 I + v$ for some $\lambda_1\in \Bbb R$ and $v\in V$.  
If $\lambda_1 = 0$ then $k=v\neq 0$ and $v\in R(h_0)T\otimes 
TR(h_0)$, which contradicts $(ii)$.  If $\lambda_1 \neq 0$ then 
$\frac{k}{\lambda_1} = I + \frac{v}{\lambda_1} \notin \{\lambda 
h_0: \lambda \in \Bbb C\}$ and $\frac{k}{\lambda_1} \in
(R(h_0)T\otimes TR(h_0))\cap \{\lambda I+V: 
\lambda \in \Bbb C\}$.  We have $\frac{k}{\lambda_1}-h_0= 
I+\frac{v}{\lambda_1} - (I-(I-h_0)) = \frac{v}{\lambda_1}+(I-
h_0)\in V$ and $\frac{v}{\lambda_1}+(I-h_0)\neq 0$; moreover
$\frac{v}{\lambda_1}+(I-h_0)\in R(h_0)T\otimes TR(h_0)\cap V$.  
This contradicts $(ii)$.

$(iii)\implies (i):$  If $h_0$ is not extremal then there are two 
marginal tracial states $h_1$ and $h_2$ on $B\otimes B$ such that 
$h_0 \neq h_i$ for $i=1,2$ and $h_0 = \frac{h_1+h_2}{2}$.  
Therefore $0 \leq h_i \leq 2h_0$ for $i=1,2$, and so by 
the preceding lemma it 
follows that $h_i\in (R(h_0)T\otimes TR(h_0))\cap \{\lambda I+V: 
\lambda \in \Bbb C\}$ for $i=1,2$.  If $h_1=\lambda h_2$ then 
$\tau(h_1) = \lambda \tau(h_2) = \lambda$, hence $\lambda =1$ and 
so dim$((R(h_0)T\otimes TR(h_0))\cap \{\lambda I+V: 
\lambda \in \Bbb C\}) \geq 2$.\qed
\enddemo
\remark{Remark 2.2}  In \cite{P1} K. R. Parthasarathy provides
necessary and sufficient 
conditions for a state, with restrictions 
$\rho_1$ on $B\otimes I$ and $\rho_2$ on $I\otimes B$, to
be extremal among all states with the same restrictions.  In
the 
appendix to this paper we show that our condition $(ii)$ above
is equivalent to
Parthasarathy's condition for marginal tracial states.
\endremark

\proclaim{Corollary 2.4} (cf. \cite{P1}) If $h_0$ is an 
extremal marginal tracial
state then $R(h_0) < I$, i.e., $h_0$ is not 
invertible.
\endproclaim
\demo{Proof}  If $R(h_0)=I$ then $(R(h_0)T\otimes TR(h_0))\cap V 
= V = \{0\}.$  On the other hand, $V= N\otimes N$ where $N=\{a\in 
B: \tau_B(a)=0\}$, where $\tau_B$ is the tracial state on $B$.  
Hence $V\neq \{0\}$.\qed
\enddemo

Let $\Cal H$ be the set of all Hilbert-Schmidt class matrices of 
$B$; then $\Cal H = B$ as sets and the inner product of $\Cal H$ is 
given by $\left< a,b\right> = 
\tau_B(b^*a)$, for $a,b\in B$, where $\tau_B$ is the tracial 
state on $B$.  The norm on $\Cal H$ is given by 
$\Vert a \Vert_2 = 
\tau_B(a^*a)^{\frac12}, a\in \Cal H$.   
$\Cal H\otimes \Cal H$ coincides with $B\otimes 
B$ as a set and $\Cal H\otimes \Cal H$ is a Hilbert space
with inner product given by 
$\left< a\otimes b,c\otimes d \right>=\tau( (c\otimes 
d)^*(a\otimes b))$ for $a,b,c,d\in B$ and extended to all of
$\Cal H \otimes \Cal H$ by linearity.
Moreover, since $B= \Bbb CI + N$ it follows that
$\Cal H\otimes \Cal H = B\otimes B = \Bbb CI \oplus 
(N\otimes I + I\otimes N) \oplus (N\otimes N)$ where we abuse 
notation slightly by using $I$ to 
denote the both identity on $B$ and on $B\otimes B$.

\proclaim{Lemma 2.5}  
$ \{\lambda I: \lambda \in \Bbb C\}\oplus 
(N\otimes I + I\otimes N) \oplus (N\otimes N)$ is an orthogonal
decomposition of $\Cal H\otimes \Cal H$.
\endproclaim
\demo{Proof}  For $\lambda \in \Bbb C$ 
we have $\tau( (\lambda I)^*(N\otimes I + I\otimes 
N))= \bar{\lambda}\tau(N\otimes I + I\otimes N) = 0$ and 
$\tau((N\otimes I+I\otimes N)^*(N\otimes N))= \tau(N^*N\otimes 
N+N\otimes N^*N)= \tau(N^*N\otimes I)\tau(I\otimes 
N)+\tau(N\otimes I)\tau(I\otimes N^*N) =0$.\qed
\enddemo

Henceforth we shall often use the notation $B$ rather 
than $T$ or $\Cal H$ because they are the same as sets.  

Let $E$ be a linear subspace of $B\otimes B$.  We shall use the 
notation $E^o$ to denote the orthogonal complement of $E$ with 
respect to the scalar product defined by the tracial state
 $\tau$ on 
$B\otimes B$.  We shall also simplify notation by using $R$ to 
denote the range projection $R(h_0)$ 
of a fixed extremal marginal tracial state $h_0$.  
Then we have the following lemma.

\proclaim{Lemma 2.6}  $B\otimes B = \Bbb CI + (N\otimes I + 
I\otimes N) + (RB\otimes BR)^o$.
\endproclaim
\demo{Proof}  By Theorem $2.3$, $(RB\otimes BR)\cap V=\{0\}$; hence 
$B\otimes B = ((RB\otimes BR)\cap V)^o = (RB\otimes BR)^o + V^o = 
(RB\otimes BR)^o+(B\otimes I + I\otimes B)$.\qed
\enddemo

The following lemma is an immediate consequence of the preceding 
lemma.

\proclaim{Lemma 2.7}  $R(B\otimes I + I\otimes B)R = \Bbb CR+ 
R(N\otimes I + I\otimes N)R = RB\otimes BR$.
\endproclaim

\remark{Remark $2.3$}Note that the conclusion of the lemma shows that
$R\neq I$ and so we recapture the result of Corollary $2.4$.  
\endremark   

The identification made in the preceding lemma between
$RB\otimes BR$ and
$\Bbb CR+R(N\otimes I + I\otimes N)R$ is trivial when
 dim($R$)$=1$ and seems to be rather puzzling otherwise, since
it does not seem intuitive to us that $\Bbb CR+R(N\otimes I + I\otimes 
N)R$ is isomorphic to a full matrix algebra for dim($R$)$>1$.  
For this reason
the lemma suggests that any extremal marginal tracial state 
must in fact be pure.

In what follows we shall appeal to the following 
observation.   
Since $h_0\in RB\otimes BR$ and $\tau(h_0R(N\otimes I + I\otimes 
N)R) = \tau((N\otimes I + I\otimes N)h_0) = 0$, we have the 
orthogonal decomposition:
$$
RB\otimes BR = \Bbb Ch_0 \oplus R(N\otimes I + I\otimes N)R.
$$
Moreover, if $W = \{a\in RB\otimes BR: \tau(a) = 0\}$ then 
by the preceding equation $W= 
h_0^{\frac12}R(N\otimes I + I\otimes N)Rh_0^{\frac12} =
h_0^{\frac12}(N\otimes I + I\otimes N)h_0^{\frac12}$.

Our next main goal, in Theorem $2.9$, is to obtain
necessary and 
sufficient conditions for a extremal marginal tracial state
 to be a pure 
state.  To prove this result it will be helpful to study 
two norms on the linear subspace $\Bbb CI + V$ of $B\otimes B$.  
First note that since $B\otimes B$ has the orthogonal 
decomposition $B\otimes B = (\Bbb CI + V)\oplus (N\otimes I + 
I\otimes N)$ we may identify $\Bbb CI + V$ with the quotient 
space $B\otimes B/(N\otimes I + I\otimes N)$ of $B\otimes B$.

The first norm we impose on $\Bbb CI + V$ is the $C^*$-norm 
$\Vert \cdot \Vert$ which 
the space inherits as a subspace of $B\otimes B$.  The other is the 
quotient norm $\Vert \vert \cdot \vert \Vert$ given by
$\Vert \vert a \vert \Vert = 
\underset{y\in N\otimes I + I\otimes N}\to{inf} \Vert a+y \Vert,
a\in \Bbb CI + V$.
Clearly $\Vert \vert a \vert \Vert \leq \Vert a \Vert$, all $a\in 
\Bbb CI + V$.  

Let $\Vert \cdot \Vert_1$ be the trace norm on $T\otimes T$.  
Then $\Vert a \Vert_1 = \tau((a^*a)^{\frac12}), a\in T\otimes T$, 
and is the dual norm of the $C^*$-norm $\Vert \cdot \Vert $ on 
$B\otimes B$.  Let $\Vert \cdot \Vert^*$ be the dual norm on $\Bbb 
CI + V (\subset T\otimes T)$ with respect to the norm $\Vert 
\cdot \Vert$ on $\Bbb CI + V (\subset B\otimes B)$.  By the 
general theory of Banach spaces 
$$
(B\otimes B/(N\otimes I + I\otimes N))^* = (N\otimes I + I\otimes 
N)^o = \Bbb CI + V \subset T\otimes T
$$
with respect to the norm $\Vert \cdot \Vert_1$; hence we have, 
for $f\in \Bbb CI + V$, 
$$
\align
\Vert f \Vert_1 &= \underset{\Vert \vert x \vert \Vert \leq 1,x\in \Bbb 
CI+V}\to{\text{sup}} \vert f(x) \vert 
 = \underset{\Vert  x \Vert \leq 1,x\in 
B\otimes B}\to{\text{sup}} \vert f(x) \vert, \qquad \text{and} \\
\Vert f \Vert^* &= \underset{\Vert x \Vert \leq 1,x\in 
\Bbb CI+V}\to{\text{sup}} \vert f(x) \vert,
\endalign
$$
and so $\Vert f \Vert^* \leq \Vert f \Vert_1$.

Suppose $f\in \Bbb CI+V$ is a state on $B\otimes B$; then
$1 = \Vert f \Vert_1 = f(I)$.  By restriction $f$ may be viewed 
as a linear functional on $\Bbb CI + V$, and we have
$\Vert f \Vert_1 = 
\underset{\Vert \vert x \vert \Vert \leq 1,x\in \Bbb 
CI+V}\to{\text{sup}} \vert f(x) \vert$.  Let $\Delta$ be the 
set of all such linear functionals.  Then $\Delta$ coincides with 
the set of all positive elements $h$ in $\Bbb CI+ V$ such that
$\tau(h) = 1$.  Since
$f(I) \leq \Vert f 
\Vert^* \leq \Vert f \Vert_1$, $\Vert f \Vert^* = \Vert f 
\Vert_1$ for $f\in \Delta$.

Let $h_0$ be an extreme point of $\Gamma$ and let $I(h_0)$ be the 
set of all states $k$ on $B\otimes B$ such that $\tau(ah_0) = 
\tau(ak)$ for $a\in \Bbb CI+V$; then $k-h_0\in (\Bbb CI + V)^o = 
N\otimes I + I \otimes N$ and therefore there is a selfadjoint 
element $a$ in $N\otimes I + I\otimes N$ such that $k = h_0 + a$.  
Conversely we have the following:

\proclaim{Lemma 2.8} Let $h$ be an arbitrary state on $B\otimes B$ and 
consider its restriction to $\Bbb CI + V$; then there is a  
unique selfadjoint element $\ell$ in $\Bbb CI + V$ such that 
$\tau(xh) = \tau(x\ell)$ for $x\in \Bbb CI + V$.  In fact
$\ell = 
( I-(P-Q)^2)(h)$.
\endproclaim
\demo{Proof}  
  A straightforward calculation shows that $(P-Q)^2(P-Q)^2=(P-
Q)^2$, and $((P-Q)^2)^* = (P-Q)^2$ so that $(P-Q)^2$ is the 
orthogonal projection of $B\otimes B$ onto $N\otimes I + I\otimes 
N$ in the Hilbert space $B\otimes B$.  Hence $I-(P-Q)^2$ is the 
orthogonal projection of $B\otimes B$ onto $\Bbb CI+V$.  Hence 
$\ell = 
( I-(P-Q)^2)(h)\in \Bbb CI+V$.  

Next we show that $\tau(xh) = \tau(x\ell)$ for all $x\in \Bbb CI+V$.
To see this note that for any $a,b\in B\otimes B$, $\tau(P(a)b) = 
\tau(P(a)P(b)) = \tau(aP(b))$.  Analogously $\tau(Q(a)b) = 
\tau(aQ(b))$.  Therefore 
$\tau ( (P-Q)^2(a)b ) = \tau ( (P-Q)(a)(P-
Q)(b)) = \tau(a(P-Q)^2(b))$, and so
$\tau (( I-(P-Q)^2)(a)b) =
\tau(a(I-(P-Q)^2)(b))$.  Now replace $a$ with $x\in \Bbb CI+V$ and
$b$ with $h$ in this equation.  From the preceding paragraph
$x=(I-(P-Q)^2)(x)$, so that $\tau(xh) = \tau(x\ell)$ for all
$x\in \Bbb CI + V$.  The uniqueness of $\ell$ is straightforward.
\qed
\enddemo

The following gives a characterization of the extreme points of 
$\Gamma$ that are in fact pure states.  As pointed out in 
\cite{P1} all extreme points are pure in the case $n=2$, i.e. 
when the algebra $B$ is isomorphic to the matrix algebra 
$M_2(\Bbb C)$.

\proclaim{Theorem 2.9}  Let $h_0$ be an extreme point of 
$\Gamma$, the set of marginal tracial states on $B\otimes B$.  
Then $h_0$ is a pure state on $B\otimes B$ if and only if 
$(I-(P-Q)^2)(RB\otimes BR) \subset RB\otimes BR$.
\endproclaim

\demo{Proof}
Suppose $h_0$ is a pure state; then $RB\otimes BR =\Bbb CR$ and 
so $h_0=n^2R$.  By Lemma $2.1$, $(P-Q)(h_0)=0$, so
$$
(I-(P-Q)^2)(R)=\frac{1}{n^2}(I-(P-Q)^2)(h_0)) = 
\frac{1}{n^2}h_0=R.
$$
Hence $(I-(P-Q)^2)\Bbb CR = \Bbb CR$.

For the converse suppose $(I-(P-Q)^2)(RB\otimes BR)\subset 
RB\otimes BR$.
By the proof of Lemma $2.8$, for any  $h\in RB\otimes BR$,
$(I-(P-Q)^2)(h)\in \Bbb CI+V$.  Therefore
$(I-(P-Q)^2)(h)\in (\Bbb CI+V)\cap (RB\otimes BR)$.  Combining 
this with the results of Theorem $2.3$ shows 
$(I-(P-Q)^2)(h)\in \{\lambda h_0: \lambda \in \Bbb C\}$, and so
$(I-(P-Q)^2)(a) \in \{\lambda h_0:\lambda\in \Bbb C\}$ for $a\in 
RB\otimes BR$; hence
$$
(I-(P-Q)^2)(a)=\tau(a)h_0,\qquad a\in RB\otimes BR. \tag{2.1}
$$
By Lemma $2.1$, $P(h_0) = I = Q(h_0)$. Therefore
$(P-Q)^2(h_0) = 0$ so that we may rewrite the equation above
as 
$$
a-\tau(a)h_0 = (P-Q)^2(a-\tau(a)h_0).  \tag{2.2}
$$ 
Hence $a-\tau(a)h_0\in N\otimes I + I\otimes 
N$, $a\in RB\otimes BR$.

Now assume (to obtain a contradiction) that
$h_0$ is not a pure state. Then if $r= 
\text{rank}(R)$, $r > 1$, $RB\otimes BR$ is isomorphic to the 
algebra of $r\times r$ matrices, and therefore the linear space
$W= \{a\in RB\otimes BR:\tau(a) = 0 \}$ has dimension $r^2-1$.  
From the preceding paragraph $W$ is also a subset of $N\otimes I 
+ I\otimes N$.  Let $a \neq 0$ be a selfadjoint element of $W$
satisfying $\Vert a \Vert_2 =1$, then 
there are selfadjoint elements $\ell, m$ of $N$ such that
$a= \ell\otimes I+I\otimes m$.  

We show $2\ell\otimes m = h_0 -I$.  To see this apply Eq. 2.2
to $a^2$ to see that 
$$
a^2-h_0 = (P-Q)^2(a^2-h_0) \tag{2.3}.
$$
But
$$
\align
(P-Q)^2(a^2 - h_0) &= (P-Q)^2(a^2) \\
&= (P-Q)^2(\ell^2\otimes I + I\otimes m^2+2\ell \otimes m)\\
&= (P-Q)^2(\ell^2\otimes I + I\otimes m^2) \\
&= (P-Q)^2(\ell^2\otimes I + I\otimes m^2 - I) \\
&= \ell^2\otimes I + I\otimes m^2 - I, \\
\endalign
$$
where the last equality follows from the observation that
$\ell^2\otimes I + I\otimes m^2 - I\in N\otimes I + I\otimes N$, 
since
\roster
\item"{(i)}"
$\tau(\ell^2\otimes I + I\otimes m^2 - I) = \tau(a^2)-1=0$, and
\item"{(ii)}" $\ell^2\otimes I + I\otimes m^2 - I$ is orthogonal to 
$N\otimes N$.  
\endroster
From the preceding calculation and Equation 2.3 we have
$2\ell \otimes m = h_0 - I$ for all self-adjoint $a$ in $W$ with 
$\Vert a \Vert_2 = 1$.  

By Corollary $2.4$, $h_0 \neq I$ hence $\ell \otimes m \neq 0$, 
hence $\ell \neq 0$ and $m \neq 0$.  Recall \cite{T} that for any 
functional $\rho$ on $B$ the map $x\otimes y \rightarrow 
\rho(y)x$, $x,y\in B$ extends by linearity to a well-defined 
(right) slice map from $B\otimes B$ to $B$. If we define $\rho$  
by $\rho(y) = \tau(m^*y), y\in B$
then the corresponding right slice 
map sends $\ell \otimes m$ to $\tau(m^*m)\ell \neq 0$.  Similarly 
there is a left slice map that sends $\ell \otimes m$ to 
$\tau(\ell^* \ell)m \neq 0$.  But if $b = \ell_1 \otimes I + 
I\otimes m_1$ is any other selfadjoint element of 
$W$ with $\Vert b \Vert_2 = 1$ then $2\ell_1\otimes m_1 = h_0-I = 
2\ell \otimes m$, so applying the left and right slice maps 
above to 
$\ell_1 \otimes m_1$ shows that $\ell_1$ (resp., $m_1$) is a 
scalar multiple of $\ell$ (resp. of $m$).   
It follows immediately that if $c= 
\ell_2 \otimes I + I\otimes m_2$ is {\it any} selfadjoint element of 
$W$, then independent of $\Vert c \Vert_2$, $\ell_2$ and $m_2$ are
scalar multiples of $\ell$ and $m$ respectively.  Since $W = W^*$ 
it then follows that for any element $d$ of $W$, $d$ lies in the 
subspace of $N\otimes I + I\otimes N$ spanned by $\ell \otimes 
I,I\otimes m$.  Hence $W$ is at most two-dimensional.  Since 
$W$ has dimension $r^2-1$, the only possibility is
$r = 1$.  Hence $h_0$ must be a pure state 
on $B\otimes B$. \qed
 \enddemo

\remark{Problem}  Let $h_0$ be an extreme point of $\Gamma$ and 
let $h$ be an arbitrary state on $B\otimes B$ with $h\in 
RB\otimes BR$.  Can we conclude that $\Vert (I-(P-
Q)^2)(h)\Vert_1 \leq 1$?.  
If so then $h_0$ is a pure state.  In fact, $\tau((I-(P-
Q)^2)(h))=\tau(h) = 1$ and so $(I-(P-Q)^2)(h)$ is a state.  Since
$(I-(P-Q)^2)(h)\in \Bbb CI+V$, it is a marginal tracial state.
  Now put $h_0 =
\sum_{j=1}^m \lambda_j n^2e_j$, where $\lambda_j > 0$ for
$j=1,2,\dots , m$ with $\sum_{j=1}^m \lambda_j = 1$ and $\{e_j: 
j=1,2,\dots , m\}$ is a family of mutually orthogonal 
one-dimensional projections 
in $B\otimes B$
such that $\sum_{j=1}^m e_j = R$.  
Then $(I-(P-Q)^2)(h_0) = h_0 = \sum_{j=1}^m \lambda_j n^2 (I-(P-
Q)^2)(e_j)$.  Since $h_0$ is extreme in $\Gamma$, $n^2(I-(P-
Q)^2)(e_j) = h_0$ for all $j$; hence $(I-(P-Q)^2)(R) =
(I-(P-Q)^2)\left( \sum_{j=1}^m e_j \right) = \frac{m}{n^2} h_0$ 
and so $h_0 = \frac{n^2}{m} R$.  Hence for an arbitrary 
projection $p$ of $RB\otimes BR$, $(I-(P-Q)^2)(p) = \tau(p)h_0$,
and so by Theorem $2.9$, $h_0$ is pure.
\endremark

\proclaim{Corollary 2.10}  Let $h_0$ be an extremal marginal 
tracial state on $B\otimes B$.  Then the following conditions are 
equivalent.
\roster
\item"{(i)}"  $h_0$ is a pure state on $B\otimes B$.
\item"{(ii)}"  for any $k\in RB\otimes BR$ the restriction of $k$ 
to $\Bbb CI + V$ is $\tau(k){h_0}_{\restriction{\Bbb CI+V}}$.
\item"{(iii)}" for $\ell\in h_0^{1/2}(N\otimes I + I \otimes 
N)h_0^{1/2}$, the restriction of $\ell$ to $\Bbb CI+V$ is $0$.
\item"{(iv)}" $h_0^{1/2}(N\otimes I + I\otimes N)h_0^{1/2} 
\subset (N\otimes I + I\otimes N)$.
\item"{(v)}" $h_0^{1/2}(\Bbb CI+V)h_0^{1/2} \subset \Bbb CI + V$.
\endroster
\endproclaim
\demo{Proof}  $(i)\implies (ii)$:  $R=R(h_0)$ is a rank one 
projection so $h_0 = n^2R$ and $k=\lambda R$. 
Then 
$\tau(k) = \lambda \tau(R) = \lambda \frac{1}{n^2}$, so $k= 
\frac{\lambda}{n^2}h_0$.

$(ii)\implies (iii)$:  By the remark following Lemma $2.7$,
since $h_0^{1/2}(N\otimes I + I\otimes 
N)h_0^{1/2} = W = \{a\in RB\otimes BR: \tau(a)=0\}$, so by $(ii)$ 
the restriction of $\ell$ to $\Bbb CI+V$ must be $0$.

$(iii)\implies (iv)$:  For $a\in \Bbb CI+V$, by Lemma $2.8$ 
we have $\tau(a(I-(P-Q)^2)(h_0^{1/2}(N\otimes I + I\otimes 
N)h_0^{1/2}) = 0$.  On the other hand, 
$\tau(a(I-(P-Q)^2)(h_0^{1/2}(N\otimes I + I\otimes 
N)h_0^{1/2})= \tau((I-(P-Q)^2)(a)(h_0^{1/2}(N\otimes I+I\otimes 
N)h_0^{1/2}) = \tau(ah_0^{1/2}(N\otimes I+I\otimes 
N)h_0^{1/2}) = 0$ for $a\in \Bbb CI+V$; hence 
$h_0^{1/2}(N\otimes I + I\otimes N)h_0^{1/2}\subset (\Bbb CI+V)^o 
= N\otimes I + I\otimes N$.

$(iv)\implies (v)$:  By the above equality $0 = \tau((\Bbb 
CI+V)h_0^{1/2}(N\otimes I + I\otimes N)h_0^{1/2}) = 
\tau(h_0^{1/2}(\Bbb CI+V)h_0^{1/2}(N\otimes I+I\otimes N))$; 
hence
$h_0^{1/2}(\Bbb CI+V)h_0^{1/2}\subset (N\otimes I + I\otimes 
N)^o=(\Bbb CI+V)$.

$(v)\implies (iv)$: clear.

$(iv)\implies (iii)$:  $W = h_0^{1/2}(N\otimes I + I\otimes 
N)h_0^{1/2} \subset (N\otimes I + I \otimes N)$; hence
$(I-(P-Q)^2)(W)\subset (I-(P-Q)^2)(N\otimes I + I\otimes N) = 0$.

$(iii)\implies (i)$:  $RB\otimes BR = \Bbb Ch_0 + W$ and 
$(I-(P-Q)^2)(\Bbb Ch_0+W) = \Bbb Ch_0 \subset RB\otimes BR$; hence 
by Theorem $2.9$, $h_0$ is a pure state.  This completes the 
proof. \qed
\enddemo

\subhead
3. Appendix
\endsubhead

Let $\rho_1$ and $\rho_2$ be states on $B$ and let $\Cal 
E(\rho_1,\rho_2)$ be the convex set consisting
of all states $\rho$
on $B\otimes B$ 
which restrict on $B\otimes I$ to $\rho_1$ and on $I\otimes B$ to 
$\rho_2$, i.e., $\rho(x\otimes I) = \rho_1(x)$ and $\rho(I\otimes 
y) = \rho_2(y)$, for $x,y\in B$.  In \cite{P1} Parthasarathy 
gives a necessary and sufficient condition for $\rho$ to be 
extremal in $\Cal E(\rho_1,\rho_2)$.  We show how 
Parthasarathy's criterion is related to ours in the case 
$\rho_1=\rho_2=\tau_B$ by showing that it is equivalent to 
condition $(ii)$ of Theorem $2.3$.  

The following lemma can be found in \cite{P2} (see also \cite{P1}).

\proclaim{Lemma A.1}  
Let $h$ be a state on 
$B\otimes B$ with rank $r < m$, where $m = n^2$.
Then there exists a positive invertible $r\times r$
matrix $K$, a permutation matrix $\sigma$ of $B\otimes B$, and
an $r\times (m-r)$ matrix A such that
$$
\sigma h \sigma^{-1} = \bmatrix 
K & KA \\
A^*K & A^*KA \\
\endbmatrix
$$                     
\endproclaim

\proclaim{Theorem A.2}  (cf. \cite{P1}) 
Let $h$ be a marginal tracial state on 
$B\otimes B$ with rank $r < m$.  Then $h$ is not an extremal 
marginal tracial state if and only if there is a selfadjoint 
$r\times r$ matrix $L$ such that 
$$
\sigma^{-1}\bmatrix 
L & LA \\
A^*L & A^*LA \\
\endbmatrix \sigma
$$
is in $V$.
\endproclaim
\demo{Proof}
By condition 
$(ii)$ of Theorem $2.3$
the state $h$ on 
$B\otimes B$ is not 
extremal if and only if $(RT\otimes TR)\cap V$ is nontrivial, 
where $R = R(h)$  and $V = N\otimes N$.  
Let $\sigma$ be the permutation matrix in $B\otimes B$ such that 
$h^{\sigma}= \sigma h \sigma^{-1}$ has the form of the matrix in 
the
lemma, and let $R^{\sigma} = \sigma R(h) \sigma^{-1}$.
It follows from 
spectral theory that 
$R^{\sigma} = \underset{j\rightarrow \infty}\to{lim} 
h_j^{\sigma}$
where $h_1^{\sigma} = h^{\sigma}, h_{j+1}^{\sigma} = 
(h_j^{\sigma})^{1/2}$.  Moreover, for any positive matrix $C$ in
$B\otimes B$ (or in fact, for 
any positive operator $C$ on a Hilbert space) we may 
obtain $C^{1/2}$ as a limit
$C^{1/2} = \underset{j\rightarrow \infty}\to{\text{lim}} C_j$ 
where $C_0 = 0$ and $C_{j+1} = C_j + \frac12 (C-C_j^2)$ 
\cite{Sz}.  Therefore $C^{1/2}$ is a limit of linear 
combinations of powers of $C$.  It follows that the successive
square roots 
$h_j^{\sigma}$ of the lemma
all have the form
$$
\bmatrix
Z_j & Z_jA\\
A^*Z_j & A^*Z_jA \\
\endbmatrix
$$
and therefore there is an $r\times r$ matrix $Q$ such that
the projection $R^{\sigma}$ has the form 
$$
\bmatrix
Q & QA\\
A^*Q & A^*QA \\
\endbmatrix.
$$
If $(RT\otimes TR)\cap V$ is nontrivial then there is a nonzero 
selfadjoint element $D$ in $(RT\otimes TR)\cap V$ and $D^{\sigma}= 
R^{\sigma}D^{\sigma}R^{\sigma}$.  From this equation and the form of 
$R^{\sigma}$ it follows that $D^{\sigma}$ has the form 
$
\bmatrix
L & LA \\
A^*L & A^*LA \\
\endbmatrix,
$
as in the statement of the lemma, with $L=L^*$, $L\neq 0$, and 
with $D = \sigma^{-1}D^{\sigma}\sigma$ in $V$.   
 
 To prove the converse suppose there is a selfadjoint $r\times 
r$
matrix $L$ such that the matrix
$D = \sigma^{-1}\bmatrix L & LA\\
      A^*L & A^*LA \\
\endbmatrix \sigma
$ is in $V$.  Since $K$ is invertible it is not difficult to show 
that $R^{\sigma}D^{\sigma}R^{\sigma} = D^{\sigma}$ (because the 
range of the matrix $D^{\sigma}$ is contained in the range of the 
matrix $h^{\sigma}$), so that $D= RDR$ and therefore
$D\in (RT\otimes TR)\cap V$.\qed
\enddemo

\enddocument